\documentclass[12pt]{article}
\usepackage{latexsym, amssymb}
\DeclareMathAlphabet\gothic{U}{euf}{m}{n}

\makeatletter
\def\eqnarray{\stepcounter{equation}\let\@currentlabel=\theequation
\global\@eqnswtrue
\tabskip\@centering\let\\=\@eqncr
$$\halign to \displaywidth\bgroup\hfil\global\@eqcnt\z@
  $\displaystyle\tabskip\z@{##}$&\global\@eqcnt\@ne
  \hfil$\displaystyle{{}##{}}$\hfil
  &\global\@eqcnt\tw@ $\displaystyle{##}$\hfil
  \tabskip\@centering&\llap{##}\tabskip\z@\cr}

\def\endeqnarray{\@@eqncr\egroup
      \global\advance\c@equation\m@ne$$\global\@ignoretrue}

\def\@yeqncr{\@ifnextchar [{\@xeqncr}{\@xeqncr[5pt]}}
\makeatother

\begin{document}
\bibliographystyle{tom}

\newtheorem{lemma}{Lemma}[section]
\newtheorem{thm}[lemma]{Theorem}
\newtheorem{cor}[lemma]{Corollary}
\newtheorem{voorb}[lemma]{Example}
\newtheorem{rem}[lemma]{Remark}
\newtheorem{rems}[lemma]{Remarks}
\newtheorem{prop}[lemma]{Proposition}
\newtheorem{stat}[lemma]{{\hspace{-5pt}}}
\newtheorem{obs}[lemma]{Observation}
\newtheorem{defin}[lemma]{Definition}

\newenvironment{remarkn}{\begin{rem} \rm}{\end{rem}}
\newenvironment{remarkns}{\begin{rems} \rm}{\end{rems}}
\newenvironment{exam}{\begin{voorb} \rm}{\end{voorb}}
\newenvironment{defn}{\begin{defin} \rm}{\end{defin}}
\newenvironment{obsn}{\begin{obs} \rm}{\end{obs}}

\newenvironment{emphit}{\begin{itemize} }{\end{itemize}}

\newcommand{\gota}{\gothic{a}}
\newcommand{\gotb}{\gothic{b}}
\newcommand{\gotc}{\gothic{c}}
\newcommand{\gote}{\gothic{e}}
\newcommand{\gotf}{\gothic{f}}
\newcommand{\gotg}{\gothic{g}}
\newcommand{\gothh}{\gothic{h}}
\newcommand{\gotk}{\gothic{k}}
\newcommand{\gotm}{\gothic{m}}
\newcommand{\gotn}{\gothic{n}}
\newcommand{\gotp}{\gothic{p}}
\newcommand{\gotq}{\gothic{q}}
\newcommand{\gotr}{\gothic{r}}
\newcommand{\gots}{\gothic{s}}
\newcommand{\gotu}{\gothic{u}}
\newcommand{\gotv}{\gothic{v}}
\newcommand{\gotw}{\gothic{w}}
\newcommand{\gotz}{\gothic{z}}
\newcommand{\gotA}{\gothic{A}}
\newcommand{\gotB}{\gothic{B}}
\newcommand{\gotG}{\gothic{G}}
\newcommand{\gotL}{\gothic{L}}
\newcommand{\gotS}{\gothic{S}}
\newcommand{\gotT}{\gothic{T}}

\newcommand{\mn}{\marginpar{\hspace{1cm}*} }
\newcommand{\mnn}{\marginpar{\hspace{1cm}**} }

\newcommand{\mnq}{\marginpar{\hspace{1cm}*???} }
\newcommand{\mnnq}{\marginpar{\hspace{1cm}**???} }

\newcounter{teller}
\renewcommand{\theteller}{\Roman{teller}}
\newenvironment{tabel}{\begin{list}%
{\rm \bf \Roman{teller}.\hfill}{\usecounter{teller} \leftmargin=1.1cm
\labelwidth=1.1cm \labelsep=0cm \parsep=0cm}
                      }{\end{list}}

\newcounter{tellerr}
\renewcommand{\thetellerr}{(\roman{tellerr})}
\newenvironment{subtabel}{\begin{list}%
{\rm  (\roman{tellerr})\hfill}{\usecounter{tellerr} \leftmargin=1.1cm
\labelwidth=1.1cm \labelsep=0cm \parsep=0cm}
                         }{\end{list}}
\newenvironment{ssubtabel}{\begin{list}%
{\rm  (\roman{tellerr})\hfill}{\usecounter{tellerr} \leftmargin=1.1cm
\labelwidth=1.1cm \labelsep=0cm \parsep=0cm \topsep=1.5mm}
                         }{\end{list}}

\newcommand{\Ni}{{\bf N}}
\newcommand{\Ri}{{\bf R}}
\newcommand{\Ci}{{\bf C}}
\newcommand{\Ti}{{\bf T}}
\newcommand{\Zi}{{\bf Z}}
\newcommand{\Fi}{{\bf F}}

\newcommand{\proof}{\mbox{\bf Proof} \hspace{5pt}} 
\newcommand{\remark}{\mbox{\bf Remark} \hspace{5pt}}
\newcommand{\ruimte}{\vskip10.0pt plus 4.0pt minus 6.0pt}

\newcommand{\simh}{{\stackrel{{\rm cap}}{\sim}}}
\newcommand{\ad}{{\mathop{\rm ad}}}
\newcommand{\Ad}{{\mathop{\rm Ad}}}
\newcommand{\Aut}{\mathop{\rm Aut}}
\newcommand{\arccot}{\mathop{\rm arccot}}
\newcommand{\capp}{{\mathop{\rm cap}}}
\newcommand{\rcapp}{{\mathop{\rm rcap}}}
\newcommand{\diam}{\mathop{\rm diam}}
\newcommand{\divv}{\mathop{\rm div}}
\newcommand{\codim}{\mathop{\rm codim}}
\newcommand{\RRe}{\mathop{\rm Re}}
\newcommand{\IIm}{\mathop{\rm Im}}
\newcommand{\Tr}{{\mathop{\rm Tr}}}
\newcommand{\Vol}{{\mathop{\rm Vol}}}
\newcommand{\card}{{\mathop{\rm card}}}
\newcommand{\supp}{\mathop{\rm supp}}
\newcommand{\sgn}{\mathop{\rm sgn}}
\newcommand{\essinf}{\mathop{\rm ess\,inf}}
\newcommand{\esssup}{\mathop{\rm ess\,sup}}
\newcommand{\Int}{\mathop{\rm Int}}
\newcommand{\Leibniz}{\mathop{\rm Leibniz}}
\newcommand{\lcm}{\mathop{\rm lcm}}
\newcommand{\loc}{{\rm loc}}

\newcommand{\mod}{\mathop{\rm mod}}
\newcommand{\spann}{\mathop{\rm span}}
\newcommand{\one}{1\hspace{-4.5pt}1}

\newcommand{\DWR}{}

\hyphenation{groups}
\hyphenation{unitary}

\newcommand{\tfrac}[2]{{\textstyle \frac{#1}{#2}}}

\newcommand{\cb}{{\cal B}}
\newcommand{\cc}{{\cal C}}
\newcommand{\cd}{{\cal D}}
\newcommand{\ce}{{\cal E}}
\newcommand{\cf}{{\cal F}}
\newcommand{\ch}{{\cal H}}
\newcommand{\ci}{{\cal I}}
\newcommand{\ck}{{\cal K}}
\newcommand{\cl}{{\cal L}}
\newcommand{\cm}{{\cal M}}
\newcommand{\cn}{{\cal N}}
\newcommand{\co}{{\cal O}}
\newcommand{\cp}{{\cal P}}
\newcommand{\cs}{{\cal S}}
\newcommand{\ct}{{\cal T}}
\newcommand{\cx}{{\cal X}}
\newcommand{\cy}{{\cal Y}}
\newcommand{\cz}{{\cal Z}}

\newcommand{\wtozp}{W^{1,2}\raisebox{10pt}[0pt][0pt]{\makebox[0pt]{\hspace{-34pt}$\scriptstyle\circ$}}}
\newlength{\hightcharacter}
\newlength{\widthcharacter}
\newcommand{\covsup}[1]{\settowidth{\widthcharacter}{$#1$}\addtolength{\widthcharacter}{-0.15em}\settoheight{\hightcharacter}{$#1$}\addtolength{\hightcharacter}{0.1ex}#1\raisebox{\hightcharacter}[0pt][0pt]{\makebox[0pt]{\hspace{-\widthcharacter}$\scriptstyle\circ$}}}
\newcommand{\cov}[1]{\settowidth{\widthcharacter}{$#1$}\addtolength{\widthcharacter}{-0.15em}\settoheight{\hightcharacter}{$#1$}\addtolength{\hightcharacter}{0.1ex}#1\raisebox{\hightcharacter}{\makebox[0pt]{\hspace{-\widthcharacter}$\scriptstyle\circ$}}}
\newcommand{\scov}[1]{\settowidth{\widthcharacter}{$#1$}\addtolength{\widthcharacter}{-0.15em}\settoheight{\hightcharacter}{$#1$}\addtolength{\hightcharacter}{0.1ex}#1\raisebox{0.7\hightcharacter}{\makebox[0pt]{\hspace{-\widthcharacter}$\scriptstyle\circ$}}}

 \thispagestyle{empty}
 
 \begin{center}
\vspace*{1.5cm}

{\Large{\bf Uniqueness of diffusion operators }}\\[3mm] 
{\Large{\bf and capacity estimates }}  \\[5mm]
\large Derek W. Robinson$^\dag$ \\[2mm]

\normalsize{May 2012}
\end{center}

\vspace{5mm}

\begin{center}
{\bf Abstract}
\end{center}

\begin{list}{}{\leftmargin=1.7cm \rightmargin=1.7cm \listparindent=15mm 
   \parsep=0pt}
   \item
Let   $\Omega$ be  a connected open subset of $\Ri^d$.
We analyze $L_1$-uniqueness of real second-order partial differential operators 
$H=-\sum^d_{k,l=1}\partial_k\,c_{kl}\,\partial_l$  and 
$K=H+\sum^d_{k=1}c_k\,\partial_k+c_0$  on $\Omega$ where  
$c_{kl}=c_{lk}\in W^{1,\infty}_{\rm loc}( \Omega), c_k\in L_{\infty,{\rm loc}}(\Omega)$, $c_0\in L_{2,{\rm loc}}(\Omega)$
and  $C(x)=(c_{kl}(x))>0$ for all $x\in\Omega$.
Boundedness properties of the coefficients are expressed indirectly in terms of the 
balls $B(r)$  associated with the Riemannian metric $C^{-1}$  and their Lebesgue measure
$|B(r)|$.

\noindent\hspace{10mm}First we establish that if  the balls $B(r)$ are bounded,  
the T\"acklind condition $\int^\infty_Rdr\,r(\log|B(r)|)^{-1}=\infty$
is satisfied for all large $R$ and $H$ is Markov unique then $H$ is $L_1$-unique.
If, in addition, $C(x)\geq \kappa\, (c^{T}\!\otimes\, c)(x)$ for some $\kappa>0$ and almost all $x\in\Omega$, $\divv c\in L_{\infty,{\rm loc}}(\Omega)$ is upper semi-bounded and $c_0$ is lower semi-bounded then  $K$ is also $L_1$-unique.

\noindent\hspace{10mm}Secondly, if the $c_{kl}$ extend continuously to functions which are locally bounded on 
$\partial\Omega$ and if the balls $B(r)$ are bounded we characterize Markov uniqueness of $H$ in terms of local capacity estimates
and boundary capacity estimates.
For example,
 $H$ is Markov unique if and only if for each  bounded subset $A$ of $\overline\Omega$  there exist
$\eta_n \in C_c^\infty(\Omega)$ satisfying  $\lim_{n\to\infty} \|\one_A\Gamma(\eta_n)\|_1 = 0$, where 
$\Gamma(\eta_n)=\sum^d_{k,l=1}c_{kl}\,(\partial_k\eta_n)\,(\partial_l\eta_n)$,
and $\lim_{n\to\infty}\|\one_A (\one_\Omega-\eta_n )\, \varphi\|_2 = 0$
for each $\varphi \in L_2(\Omega)$ or if and only if $\capp(\partial\Omega)=0$.

\end{list}

\vfill


\vspace{0.5cm}

\noindent
\begin{tabular}{@{}cl@{\hspace{10mm}}cl}
$ {}^\dag\hspace{-5mm}$&   Mathematical Sciences Institute (CMA)    &  {} &{}\\
  &Australian National University& & {}\\
&Canberra, ACT 0200, Australia && {} \\
  &derek.robinson@anu.edu.au 
 & &{}\\
\end{tabular}

\newpage
\setcounter{page}{1}

\section{Introduction}\label{S1}

 Let   $\Omega$ be  a connected open subset of $\Ri^d$ and define the second-order divergence-form operator $H$ on the domain $D(H)=C_c^\infty(\Omega)$ by
 \begin{equation}
H=-\sum^d_{k,l=1}\partial_k\,c_{kl}\,\partial_l
\label{elcap1.2}
 \end{equation}
 where the $c_{kl}=c_{lk}$  are real-valued  functions in $ W^{1,\infty}_{\rm loc}(\Omega)$, 
  and  the  matrix $C=(c_{kl})$ is strictly elliptic, i.e.\ $C(x)>0$ for all $x\in\Omega$.
It is possible that  the coefficients can  have degeneracies
as $x\to\partial\Omega$, the boundary of $\Omega$, or  as $x\to\infty$.

The operator $H$ is defined to be $L_1$-unique if it has a unique $L_1$-closed extension which generates a strongly continuous semigroup on $L_1(\Omega)$.
Alternatively, it is defined to be Markov unique if it has a unique $L_2$-closed extension which generates a submarkovian semigroup on the spaces $L_p(\Omega)$.
Markov uniqueness is a direct consequence of  $L_1$-uniqueness since distinct submarkovian extensions give distinct $L_1$-extensions. But the converse implication is not valid in general.
The converse was established in \cite{RSi5}  for bounded coefficients $c_{kl}$ and the proof was extended in \cite{RSi4}
to allow a growth of the coefficients  at infinity.
The converse can, however, fail if the coefficients grow too rapidly (see \cite{RSi4} Section~4.1).
The principal aim of the current paper is to establish the equivalence of Markov uniqueness and 
$L_1$-uniqueness of $H$ from properties of the Riemannian geometry 
 defined by the metric $C^{-1}$ which give, implicitly, optimal growth bounds on the coefficients.

Our arguments extend to   non-symmetric operators
 \begin{equation}
 K=H+\sum^d_{k=1}c_k\,\partial_k+c_0
 \label{elcap1.1}
 \end{equation}
with the real-valued lower-order coefficients satisfying the following three conditions:
 \begin{equation}
 \left.
\begin{array}{rl}
1.& \,c_0\in L_{2,\rm loc}(\Omega) \mbox{  is lower semi-bounded, }\\[5pt]
2.& \,c_k\in L_{\infty,\rm loc}(\Omega) \mbox{ for each } k=1,\ldots,d, \:\divv c\in L_{\infty,\rm loc}(\Omega)\hspace{1cm}\\[5pt]
{}&\mbox{and } \divv c \mbox{ is upper semi-bounded, }\\[5pt]
3.& \mbox{ there is a } \kappa>0 \mbox{ such that }
C(x)\geq \kappa\, \big(c^T\!\otimes c\big)(x) \mbox{ for }\\[5pt]
{}&\mbox{  almost all } x\in\Omega.
\end{array}
\right\}\label{elcap1.12}
\end{equation}
In the second  condition  $c=(c_1,\ldots,c_d)$ and  $\divv c=\sum^d_{k=1}\partial_k c_k$ with the partial derivatives understood in the distributional sense.
The third condition  in (\ref{elcap1.12})   
is understood in the sense of matrix ordering, i.e.\ $(c_{kl}(x))\geq\kappa\,(c_k(x)c_l(x))$ for almost all $x\in \Omega$.
These conditions together with the general theory of accretive sectorial forms are sufficient to ensure that $K$  has an extension which generates a strongly continuous semigroup on $L_1(\Omega)$ (see Section~\ref{S2}).
As in the symmetric case $K$ is defined to be $L_1$-unique if it has a unique such extension.

\smallskip

The  Riemannian 
distance 
 $d(\,\cdot\,;\,\cdot\,)$ corresponding to the metric $C^{-1}$ can be defined  in various equivalent ways but in particular by 
\begin{equation}
d(x\,;y)=\sup\{\psi(x)-\psi(y): \psi\in W^{1,\infty}_{\rm loc}(\Omega)\,,\,\Gamma(\psi)\leq1\}
\label{elcap1.30}
\end{equation}
for all $x,y\in \Omega$  where $\Gamma$,  the {\it carr{\'e} du champ} of $H$, denotes the positive map 
\begin{equation}
\varphi\in   W^{1,2}_{\rm loc}(\Omega) \mapsto \Gamma(\varphi)=\sum^d_{k,l=1}c_{kl}(\partial_k\varphi)(\partial_l\varphi)
\in L_{1,{\rm loc}}(\Omega)
\;.
\label{elcap1.3}
\end{equation}
Since $\Omega$ is connected and $C>0$ it follows that $d(x\,;y)$ is finite for all $x,y\in\Omega$ but one can have $d(x\,;y)\to\infty$ as $x$, or $y$, tends to the boundary $\partial\Omega$.
Throughout the sequel we choose coordinates such that $0\in\Omega$ and 
 denote the Riemannian distance to the origin by~$\rho$.
 Thus $\rho(x)=d(x\,;0)$ for all $x\in \Omega$.
The Riemannian ball of radius $r>0$ centred at $0$ is then defined by $B(r)=\{x\in\Omega: \rho(x)<r\}$ and its volume
(Lebesgue measure) is denoted by $|B(r)|$.

There are two properties of the balls $B(r)$ which are important in our  analysis.
First, the balls $B(r)$ must be bounded  for all $r>0$.
It follows straightforwardly that this is equivalent to the condition that $\rho(x)\to\infty$ as $x\to\infty$, i.e.\ as $x$ leaves any compact subset of~$\Omega$.
Secondly, it is essential to have control of the growth of  the volume $|B(r)|$ (Lebesgue measure) of the balls.
Our results are based on the T\"acklind condition \cite{Tac},   
\begin{equation}
\int^\infty_R dr\,r(\log |B(r)|)^{-1}=\infty
\label{elcap1.31}
 \end{equation}
for all large $R$.
In particular this condition is satisfied if there are $a, b>0$ such that $|B(r)|\leq a\,e^{b\,r^2\log(1+r)}$ for all $r>0$.

T\"acklind established the Cauchy equation on $\Ri^d$  has a unique solution
 within the class of functions satisfying
a growth condition of the type (\ref{elcap1.31}).
Moreover, uniqueness can fail if the growth bound is not satisfied.
Subsequently Grigor'yan  (see \cite{Gri6}, Theorem~1, or \cite{Gri5}, Theorem~9.1) used  condition (\ref{elcap1.31}) to prove that the heat semigroup generated by the Laplace-Beltrami operator  on a geodesically complete manifold is stochastically complete, i.e.\ it conserves probability. 
But stochastic completeness of the heat semigroup  is equivalent to $L_1$-uniqueness of the Laplace-Beltrami operator (see, for example, \cite{Dav14} Section~2).
Thus (\ref{elcap1.31}) suffices for  $L_1$-uniqueness of the Laplace-Beltrami operator.
Our aim is to prove that the T\"acklind condition and a variation of  Grigor'yan's arguments  are sufficient to establish $L_1$-uniqueness of $H$ and $K$.
In our analysis   Markov uniqueness of $H$  plays the same role as geodesic completeness of the manifold.

\begin{thm}\label{ttdep1.1}
Adopt the foregoing assumptions.
Assume   the Riemannian balls $B(r)$ are bounded for all $r>0$ and the T\"acklind condition  $(\ref{elcap1.31})$ is satisfied.
Further assume that $H$ is Markov unique.
Then $H$ and $K$ are $L_1$-unique.
\end{thm}

The theorem  extends results obtained in collaboration with El Maati Ouhabaz \cite{OuR} based on  conservation arguments
which  place more restrictive restrictions on the lower-order coefficients.

Theorem \ref{ttdep1.1}
 will be proved in Section~\ref{S3} after the discussion of some preparatory material in Section~\ref{S2}.
Finally, in Section~\ref{S5} we discuss the characterization of Markov uniqueness of $H$ in terms of capacity estimates.
These latter estimates give a practical method of establishing the Markov uniqueness property.
They also establish that if the  coefficients $c_{kl}$ extend by continuity to locally bounded functions on $\overline\Omega$ then Markov uniqueness is equivalent to the capacity of the boundary of $\Omega$ being zero.

For  background information and related results on uniqueness properties of diffusion operators we refer to  
Section~3.3 of \cite{FOT} together with  the lecture notes of Eberle \cite{Ebe} and references therein. 

\section{Preliminaries }
\label{S2}

In this section we first recall some basic results on  Markov uniqueness of the symmetric operator $H$ defined by (\ref{elcap1.2}).
These results do not require any restrictions on the growth of the coefficients of $H$ or on the Riemannian geometry.
Secondly, we  discuss the accretivity properties, etc.\ of the non-symmetric operator $K$ and its Friedrichs extension
together with continuity and quasi-accretivity properties of the associated positive semigroup.
Although these results are formulated for the operators $H$ and $K$ they are to a large extent general properties of Dirichlet forms,
symmetric \cite{BH} \cite{FOT} or non-symmetric \cite{MR}.
Thirdly, we establish some basic regularity properties for solutions of the Cauchy equations associated with $H$ and $K$.
 
 \subsection{Markov uniqueness}\label{S2.1}
The   operator $H$    is positive(-definite) and symmetric on $L_2(\Omega)$.  
The corresponding positive, symmetric,  quadratic form $h$  is  given by  $ D(h)=C_c^\infty(\Omega)$ and 
\[
h(\varphi)
= \sum^d_{k,l=1}(\partial_k\varphi, c_{kl} \, \partial_l\varphi)
\]
where $(\,\cdot\,,\,\cdot\,)$ denotes the $L_2$-scalar product.
The form  is closable and its closure $h_D=\overline{h}$ determines   a 
positive self-adjoint extension, the Friedrichs' extension, $H_D$ of $H$ (see, for example,  \cite{Kat1}, Chapter~VI).
We use the notation $H_D$ since this extension  corresponds to Dirichlet  conditions on the boundary $\partial\Omega$.
The closure $h_D$ is a Dirichlet form and consequently the $H_D$ generates a submarkovian semigroup $S$. 
(For details on Dirichlet forms and submarkovian semigroups see  \cite{BH} \cite{FOT}  \cite{MR}.)
In particular $S$ extends from $L_2(\Omega)\cap L_1(\Omega)$ to a positive contraction semigroup $S^{(1)}$ 
on $L_1(\Omega)$ and the generator $H_1$ of $S^{(1)}$  is an extension of $H$.
Therefore $H$ has both a submarkovian extension and  an $L_1$-generator extension.

Next we define a second  Dirichlet form extension $h_N$ of $h$ as follows.
First the domain $D(h_N)$ of $h_N$  is specified by
\[
D(h_N)=\{\varphi\in W^{1,2}_{\rm loc}(\Omega):\Gamma(\varphi)+\varphi^2\in L_1(\Omega)\}
\]
where $\Gamma$ denotes the positive map defined by (\ref{elcap1.3}).
Then $h_N$ is given by
\[
h_N(\varphi)=\int_\Omega\Gamma(\varphi)=\|\Gamma(\varphi)\|_1
\]
for all $\varphi\in D(h_N)$.
The form $h_N$ is closed as a direct consequence of the strict ellipticity assumption $C>0$ (see \cite{RSi4}, Section~1, or \cite{OuR}, Proposition~2.1). 
The self-adjoint  operator $H_N$ associated with  $h_N$ is a submarkovian extension of $H$ which can be considered to correspond to Neumann boundary conditions.
In general the two submarkovian extensions $H_D$ and  $H_N$ of $H$ are distinct.
The significance of the forms $h_D$ and $h_N$ is that they are the minimal and maximal Dirichlet form extensions of $h$.

\begin{prop}\label{pcap2.1}
Let $k$ be a Dirichlet form extension of $h$.
Then $h_D\subseteq k\subseteq h_N$.
Thus if $K$ is the submarkovian extension of $H$ corresponding to $k$ one has $H_N\leq K\leq H_D$.

In particular, $H$ is Markov unique if and only if $h_D=h_N$.
\end{prop}
\proof\ The proposition  follows from elliptic regularity and some standard results in the theory of Dirichlet forms.
We briefly describe the proof of \cite{RSi4} which demonstrates that it is a local result
(see also \cite{FOT} Section~3.3.3, \cite{Ebe} Section~3c).

First one clearly has $h_D\subseteq k$. 
Hence $K\leq H_D$.
Secondly, since $C$ is strictly elliptic  $H$ is locally strongly elliptic. 
Then, by   elliptic regularity,  $C_c^\infty(\Omega)D(K)\subseteq D(\overline H)$ where $\overline H$ is the $L_2$-closure of $H$ (see \cite{RSi5}, Corollary~2.3, and \cite{RSi4}, Lemma~2.2).
Thirdly for each $\chi\in C_c^\infty(\Omega)$ with $0\leq \chi\leq1$ define the truncated form $k_\chi$ by
$D(k_\chi)=D(k)\cap L_\infty(\Omega)$ and $k_\chi(\varphi)=k(\varphi,\chi\varphi)-2^{-1}k(\chi,\varphi^2)$.
Then $0\leq k_\chi(\varphi)\leq k(\varphi)$ (see \cite{BH}, Proposition~4.1.1).
Moreover, if $\varphi\in D(K)\cap L_\infty(\Omega)$ then $\chi\varphi\in D(\overline H)$ and 
\[
k_\chi(\varphi)=(\varphi, \overline H \chi\varphi)-2^{-1}(H\chi,\varphi^2)
\;.
\]
But if $\chi_1\in C_c^\infty(\Omega)$ with $\chi_1=1$ on $\supp\chi$ then $\varphi_1=\chi_1\varphi\in D(\overline H)\subseteq W^{2,2}_{\rm loc}(\Omega)$, where the last inclusion again uses elliptic regularity, and
\[
k_\chi(\varphi)=(\varphi_1, \overline H \chi\varphi_1)-2^{-1}(H\chi,\varphi_1^2)=\int_\Omega \chi\,\Gamma(\varphi_1)
\]
by direct calculation.
Combining these observations one has
\[
\int_\Omega \chi\,\Gamma(\varphi_1)=k_\chi(\varphi)\leq k(\varphi)
\]
for all $\varphi\in D(K)\cap L_\infty(\Omega)$.
Then if $V$ is a relatively compact subset of $\Omega$ there is a $\mu_V>0$ such that $C(x)\geq \mu_V I$ for all $x\in V$.
Therefore choosing $\chi$ such that $\chi=1$ on $V$ one deduces that $\mu_V\int_V|\nabla\varphi|^2\leq k(\varphi)$
for each choice of $V$.
Thus $\varphi\in W^{1,2}_{\rm loc}(\Omega)$.
Moreover, $\int_V\Gamma(\varphi)\leq k(\varphi)$ for each $V$ so $\varphi\in D(h_N)$.
Consequently $D(K)\cap L_\infty(\Omega)\subseteq D(h_N)$ and 
\[
h_N(\varphi)=\sup_V\int_V\Gamma(\varphi)\leq k(\varphi)
\]
for all $\varphi\in D(K)\cap L_\infty(\Omega)$.
But since $K$ is the generator of  a submarkovian semigroup $D(K)\cap L_\infty(\Omega)$ is a core of $K$.
In addition  $D(K)$ is a core of $k$.
Therefore the last inequality extends   by continuity to all $\varphi\in D(k)$.
In particular $D(k)\subseteq D(h_N)$.
Hence $k\subseteq h_N$ and $H_N\leq K$.
\hfill$\Box$

\bigskip

The  identity $h_D=h_N$, in one guise or another, has been the basis of much of the analysis of Markov uniqueness (see, for example,
\cite{FOT}, Section~3.3, or  \cite{Ebe}, Chapter~3).
Since $h_N$ is an extension of $h_D$ the identity   is equivalent to the condition $D(h_D)=D(h_N)$.
But $D(h_D)$ is the closure of  $C_c^\infty(\Omega)$ with respect to the graph norm $\varphi\mapsto \|\varphi\|_{D(h_D)}
=(h_D(\varphi)+\|\varphi\|_2^2)^{1/2}$.
Therefore $h_D=h_N$ if and only if 
$C_c^\infty(\Omega)$ is a core of $h_N$.
Equivalently, $h_D=h_N$ if and only if $(D(h_D)\cap L_\infty(\Omega))_c$, the space of bounded functions in $D(h_D)$ with compact support in $\Omega$, is a core of $h_N$.

It follows from the Dirichlet form structure that the subspace $D(h_N)\cap L_\infty(\Omega)$ of bounded functions in $D(h_N)$ is an algebra and a core of $h_N$.
Similarly $D(h_D)\cap L_\infty(\Omega)$ is an algebra and a core of $h_D$.
The following observation on the algebraic  structure is useful for various estimates. 

\begin{prop}\label{pcap2.10}
The subalgebra $D(h_D)\cap L_\infty(\Omega)$  of $D(h_N)\cap L_\infty(\Omega)$ is an  ideal, i.e.\
\[
(D(h_D)\cap L_\infty(\Omega))\,(D(h_N)\cap L_\infty(\Omega))\subseteq D(h_D)\cap L_\infty(\Omega)\;.
\]
\end{prop}
\proof\
If $\eta\in D(h_D)\cap L_\infty(\Omega)$ then there is a sequence $\eta_n\in C_c^\infty(\Omega)$, 
with $\|\eta_n\|_2\leq \|\eta\|_2$, which converges to $\eta$ in the $D(h_D)$-graph norm.
But   if $\varphi\in D(h_N)\cap L_\infty(\Omega)$ then $\eta_n\,\varphi\in W^{1,2}_0(\Omega)$.
Further
\[
\lim_{n\to\infty}\|\eta_n\,\varphi-\eta\,\varphi\|_2\leq \lim_{n\to\infty}\|\eta_n-\eta\|_2\|\varphi\|_\infty=0
\;.
\]
Moreover,
\[
h_D(\eta_n\varphi-\eta_m\varphi)\leq 2\,h_D(\eta_n-\eta_m)\,\|\varphi\|_\infty^2+2\int_\Omega\Gamma(\varphi)(\eta_n-\eta_m)^2
\;.
\]
Since $\Gamma(\varphi)\in L_1(\Omega)$ and $\eta_n$ is $L_2$-convergent it follows by equicontinuity that
$\eta_n\,\varphi$ converges to $\eta\,\varphi$ in the $D(h_D)$-graph norm.
Thus $\eta\,\varphi\in D(h_D)$.
\hfill$\Box$

\bigskip

Although   $D(h_N)\cap L_\infty(\Omega)$  is a core of $h_N$ it does not  follow without further assumptions that 
$(D(h_N)\cap L_\infty(\Omega))_c$, the subspace   of functions  with compact support in $\overline\Omega$,  is  a core of $h_N$.
Maz'ya gives an example with $\Omega=\Ri^d$ for which this property fails (see, \cite{Maz}, Theorem~3 in Section~2.7).
We will return to the discussion of this topic in Section~\ref{S5}.

\subsection{Accretivity and continuity properties}\label{S2.2}
 
  Next we consider the non-symmetric operator $K$ defined by (\ref{elcap1.1}) with the lower order coefficients satisfying the three conditions of (\ref{elcap1.12}). 
In this subsection $K$ is viewed as an operator on the space of complex $L_2$-functions.
 Our aim is to establish accretivity and sectorial estimates  which suffice  to deduce that $K$ has a Friedrichs' extension
 which generates a strongly continuous semigroup $T$ on $L_2(\Omega)$ and that the semigroup extends to the corresponding $L_p$-spaces.
 These estimates apply equally well to the formal adjoint $K^\dagger$ of $K$.
The latter operator is defined as  the restriction of the $L_2$-adjoint $K^*$ of $K$ to $C_c^\infty(\Omega)$.
Therefore $K^\dagger$ is  obtained from $K$ by the replacements
$c\to-c$ and $c_0\to c_0-\divv c$.

After deriving the accretivity estimates we derive a local strong continuity property  for the semigroup $T$ 
 and  the dual group $T^*$ generated by the Friedrichs' extension of $K^\dagger$ both acting on $L_\infty(\Omega)$.

First define 
$L$ and $M$ on $C_c^\infty(\Omega)$ by 
\[
L\varphi=\sum^d_{k=1}c_k\partial_k\varphi\;\;\;\;\;\;{\rm and}\;\;\;\;\;\;\;
M\varphi=c_0\varphi\;.
\]
Then $K=H+L+M$.
 Let $k$ denote the corresponding sesquilinear form and quadratic form, i.e.\ $D(k)=C_c^\infty(\Omega)$,  $k(\varphi,\psi)=(\varphi, K\psi)$ and $k(\varphi)=k(\varphi,\varphi)$
 for $\varphi,\psi\in D(k)$.
 Further let $k^*$ denote the adjoint form, i.e.\ $D(k^*)=D(k)$ and $k^*(\varphi, \psi)=k(\psi,\varphi)$.
 The real  part and imaginary parts of $k$  are defined by $\Re k=2^{-1}(k+k^*)$ and $\Im k=(2i)^{-1}(k-k^*)$,
 respectively.
 In particular 
 \begin{eqnarray}
 (\Re k)(\varphi)&=&h(\varphi)+(\varphi, (c_0-2^{-1}\divv c)\varphi)\nonumber\\[5pt]
 &\geq&h(\varphi)+(\omega_0-2^{-1}\omega_1)\|\varphi\|_2^2
 \label{eacc1}
 \end{eqnarray}
 for all $\varphi\in C_c^\infty(\Omega)$ where  $ \omega_0=\essinf _{x\in\Omega}\;c_0(x)$ and $\omega_1=\esssup _{x\in\Omega}\;(\divv c)(x)$.
 Thus $\Re k$ is the form of a lower semi-bounded symmetric operator  and consequently closable.
Moreover, if $\omega= (\omega_0-2^{-1}\omega_1)$ then $k+\sigma$ is an accretive form for all $\sigma\geq -\omega$.
Next 
 \[
 (\Im k)(\varphi)=(2i)^{-1}\Big((\varphi,L\varphi)-(L\varphi,\varphi)\Big)
 \]
 for all $\varphi\in C_c^\infty(\Omega)$.
 Hence
 \begin{eqnarray} 
 | (\Im k)(\varphi)|&\leq& \|\varphi\|_2\,\|L\varphi\|_2\leq \kappa^{-1/2}\|\varphi\|_2\,h(\varphi)^{1/2}\nonumber\\[5pt]
 &\leq& (\varepsilon\,h(\varphi)+(4\varepsilon\kappa)^{-1}\|\varphi\|_2) 
  \label{eacc2}
\end{eqnarray} 
for all $\varphi\in C_c^\infty(\Omega)$ and $\varepsilon>0$
where the second step uses the third condition of (\ref{elcap1.12}).
It  follows from (\ref{eacc1}) and (\ref{eacc2}) that  $k+\sigma$ is a sectorial form for all $\sigma\geq (4\kappa)^{-1}-\omega$.
Since $\Re k$ is closable it follows that $k+\sigma$ is closable 
 with respect to the norm 
$\varphi\in C_c^\infty(\Omega)\mapsto\|\varphi\|_k=((\Re k)(\varphi)+\sigma\|\varphi\|_2^2)^{1/2}$ for any
$\sigma>-\omega$.
The closure of the form then  determines a closed extension of $K+\sigma I$ (see \cite{Kat1}, Chapter~VI or \cite{Ouh5}, Chapter~1).
Therefore by subtracting $\sigma I$ one obtains a closed extension  $K_D$ of $K$, the Friedrichs' extension.
The extension  generates a strongly continuous semigroup $T$ on $L_2(\Omega)$ which satisfies the  quasi-contractive bounds $\|T_t\|_{2\to2}\leq e^{-\omega t}$, for all $t>0$.
The estimates (\ref{eacc1}) and (\ref{eacc2})  are also valid for the adjoint form $k^*$ which is associated with the formal adjoint $K^\dagger$ of $K$.
Therefore $K^\dagger$ has a Friedrichs' extension, $K^\dagger_D=(K_D)^*$
and $K^\dagger_D$ generates the adjoint semigroup $T^*$ on~$L_2(\Omega)$.

It follows from the foregoing accretivity and sectorial properties that if $\sigma>(4\kappa)^{-1}-\omega$ then  $k+\sigma$ 
satisfies the weak sector condition I\,(2.3) of Ma and R\"ockner \cite{MR}
(see \cite{Ouh5}, Proposition~1.8).
Therefore $k+\sigma$  is accretive, closable and satisfies the weak sector condition for all sufficiently large $\sigma$.
Then   it follows from \cite{MR},  Section~II.2d, that $k+\sigma$ is a (non-symmetric) Dirichlet form.
Therefore $T$ is positive.
Moreover, $T$ extends from $L_2(\Omega)\cap L_1(\Omega)$ to a strongly continuous semigroup 
on $L_1(\Omega)$, 
and from $L_2(\Omega)\cap L_\infty(\Omega)$ to a weakly$^*$ continuous semigroup 
on $L_\infty(\Omega)$. 
Similar conclusions are valid for the adjoint form $k^*$ and the adjoint semigroup $T^*$.
Since one readily establishes that $K-(\omega_0-\omega_1)$ and $K^\dagger -\omega_0$ are both $L_1$-dissipative it then follows that 
 $\|T_t\|_{1\to1}\leq e^{-(\omega_0-\omega_1) t}$ 
and  $\|T_t\|_{\infty\to\infty}=\|T^*_t\|_{1\to1}\leq e^{-\omega_0 t}$ for all $t>0$.

One can also define an extension $K_N$ of $K$ analogous to the extension $H_N$ of $H$ by form techniques.
To this end one uses the lower semi-boundedness of $c_0$ and the third property of (\ref{elcap1.12}).
The latter ensures that the first-order operator $L$ extends to $D(h_N)$ and that the corresponding form $l$ is relatively
bounded by $h_N$ with relative bound zero.
We omit the details.

\smallskip

The weak$^*$-continuity of  the semigroup $T^{(\infty)}$ generated by $K_D$  on $L_\infty(\Omega)$ can be strengthened
by general arguments which apply equally well to the semigroup generated by $K_N$.

\begin{prop}\label{pacc1}
The  semigroup $T^{(\infty)}$ is $L_{p,\rm loc}$-continuous for all $p\in[1,\infty\rangle$.
\end{prop}
\proof\
First we prove that $T^{(\infty)}$ is $L_{1,\rm loc}$-continuous.
It clearly suffices to prove that 
\[
\lim_{t\to0}\|\one_V(I- T^{(\infty)}_t)\psi\|_1=0
\]
for all relatively compact  $V\subset \Omega$ and all positive $\psi\in L_\infty(\Omega)$.

 Let $W$ be a second relatively compact subset of $\Omega$ with $\overline V\subset W$.
  Then 
 \[
 \|\one_V(I- T^{(\infty)}_t)\psi\|_1\leq \|\one_V(I- T^{(\infty)}_t)\one_W\psi\|_1+\|\one_VT^{(\infty)}_t(\one_\Omega-\one_W)\psi\|_1
 \]
 because $\one_V(\one_\Omega-\one_W)=\one_V\one_{W^{\rm c}}=0$.
 But $\one_W\psi\in L_1(\Omega)\cap L_\infty(\Omega)$ and consequently
 \[
 \limsup_{t\to0}\|\one_V(I- T^{(\infty)}_t)\one_W\psi\|_1= \limsup_{t\to0}\|\one_V(I- T^{(1)}_t)\one_W\psi\|_1
 \leq \lim_{t\to0}\|(I- T^{(1)}_t)\one_W\psi\|_1=0
\]
by the strong continuity of $T^{(1)}$ on $L_1(\Omega)$.
Next note that $\one_V\in L_1(\Omega)$ and $\one_{W^{\rm c}}\psi\in L_\infty(\Omega)$.
But  $\one_{W^{\rm c}}\psi=(\one_\Omega-\one_W)\psi\geq0$, since $\psi\geq0$ by assumption. 
Moreover $T$ is positive.
Therefore
 \[
 \limsup_{t\to0} \|\one_VT^{(\infty)}_t(\one_\Omega-\one_W)\psi\|_1
 = \limsup_{t\to0}\,(\one_V,T^{(\infty)}_t\one_{W^{\rm c}}\psi)=0
 \]
by the weak$^*$ continuity of $T^{(\infty)}$.
Combination of these conclusions completes the proof for $p=1$.

Finally the continuity for $p\in \langle1,\infty\rangle$ follows since
\[
\|\one_V(I- T^{(\infty)}_t)\psi\|_p\leq \|\one_V(I- T^{(\infty)}_t)\psi\|_1^{1/p}((1+e^{-\omega_0t})\,\|\psi\|_\infty)^{1-1/p}
\]
by the H\"older inequality and the bounds $\|T_t^{(\infty)}\|_{\infty\to\infty}\leq e^{-\omega_0 t}$.
\hfill$\Box$

\begin{remarkn}\label{racc1} The adjoint semigroup $T^*$ is also $L_{p,\rm loc}$-continuous because it is the semigroup generated by the Friedrichs' extension $K^\dagger_D$ of the formal  adjoint $K^\dagger$ of $K$.
\end{remarkn}

\subsection{Parabolic regularity}\label{S2.3}

Next we discuss some basic regularity properties of uniformly bounded solutions of the Cauchy equations corresponding to $H$ and $K$.
The Cauchy equation is formally given by 
\[
\partial_t\psi_t+H\psi_t=0
\]
where $t>0\mapsto \psi_t$ is a function over $\Omega$ whose initial value $\psi_0$ is specified.
A precise definition will be given in the following section. 
Analysis of the Cauchy equation requires consideration of functions over the $(d+1)$-dimensional set
 $\Omega_+=\Ri_+\times\Omega$.
 We use the notation $u, v$, etc.\ for functions over $\Omega_+$ to avoid confusion with the functions $\varphi, \psi$, etc.
 over $\Omega$.
 We nevertheless use $(\,\cdot\,,\,\cdot\,)$ and $\|\cdot\|_2$ to denote the scalar product and norm on $L_2(\Omega_+)$ since this should not cause confusion.
 In particular  
 \[
 \|u\|_2=\Big(\int^\infty_0dx_0\int_\Omega dx\,|u(x_0,x)|^2\Big)^{1/2}
 \;.
 \]
 The tensor product structure ensures that the operators $H$ and $K$ and their various generator extensions act in a natural manner on  $L_2(\Omega_+)$, e.g.\ $H_D$ on $L_2(\Omega)$ is replaced by $\one_{\Ri_+}\otimes H_D$ on $L_2(\Omega_+)$.
 To avoid inessential complications we will use the same notation for the operators on the enlarged spaces, i.e.\ we identify $H_D$ with $\one_{\Ri_+}\otimes H_D$ etc.
 
We now consider the operator $ \ch=-\partial_0+H$
 acting on $C_c^\infty(\Omega_+)$.
 The formal adjoint 
  is then given by $\ch^\dagger=\partial_0+H$.
 Next we  introduce the Sobolev space 
 \[
 V^{1,2}(\Omega_+)=\{\psi\in L_2(\Omega_+): \partial_k\psi\in L_2(\Omega_+) \mbox{  for all  } k=1,\ldots d\}
 =L_2(\Ri_+)\otimes W^{1,2}(\Omega)
 \]
 and the weighted, or anisotropic, space
  \[
 V^{2,2}(\Omega_+)=\{\psi\in V^{1,2}(\Omega_+): \partial_0\psi, \partial_k\partial_l\psi\in L_2(\Omega_+) \mbox{  for all  } k,l=1,\ldots d\}
 \]
with the usual norms.
Then the spaces $ V^{-1,2}(\Omega_+)$ and $V^{-2,2}(\Omega_+)$ of distributions are defined by duality
(see, for example, \cite{Gri7} Section~6.4).
 
 The principal regularity property  used in the subsequent discussion of $L_1$-uniqueness of $H$ is the following.
 \begin{prop}\label{preg2.1}
If  $\ch^*$  denotes the $L_2$-adjoint of $\ch$ 
then $D(\ch^*)\subseteq V^{2,2}_{\rm loc}(\Omega_+)$.
 \end{prop}
\proof\
The  proposition  is a corollary of Lemma~6.19 in \cite{Gri7}.
The discussion of parabolic  regularity properties in the latter reference is for a strongly elliptic symmetric
operator $\cp$  with smooth coefficients interpreted as acting on distributions from $\cd'(\Omega_+)$. 
But since  the estimates are local only local strong ellipticity is necessary and this follows from the strict ellipticity
of the matrix $C$ of coefficients of $H$.
Moreover, the proof of  Lemma~6.19 only uses the 
assumption that the  coefficients of $\cp$ are locally Lipschitz.
Therefore the proof of Lemma~6.19 is applicable with $\cp$ replaced by $\ch^*$.
\hfill$\Box$
 
 \bigskip

In the discussion of $L_1$-uniqueness of $K$ it is convenient to introduce the operator $K_0=H+L$ on $C_c^\infty(\Omega)$
and the corresponding operator $\ck_0=-\partial_0+K_0$ on $C_c^\infty(\Omega_+)$.
Note that the formal adjoint of $K_0$ is given by $K_0^\dagger=H- L+M_0$ where $M_0$ is the operator of multiplication by the locally bounded function $-\divv c$.

 \begin{prop}\label{preg2.2}
If $\ck_0^*$  denotes the $L_2$-adjoint of $\ck_0$ 
then $D(\ck_0^*)\subseteq V^{2,2}_{\rm loc}(\Omega_+)$.
\end{prop}
\proof\
The proof of the proposition is a repetition of the argument used to prove  Lemma~6.19 in \cite{Gri7}.
The operator  $\cp$ in the latter reference is now replaced by $\ck_0^*$.
Therefore one has the terms corresponding to $\ch^*$ together with additional first-order and zero-order terms.
The additional first-order terms $-\sum^d_{k=1}c_k\partial_k$ cause no problem since they combine with the terms $-\sum_{k,l=1}^d(\partial_lc_{lk})\partial_k$. 
The zero-order term, i.e.\ multiplication by $-\divv c$, also causes no problem since $\divv c\in L_{\infty,\rm loc}(\Omega)$ by assumption.
\hfill$\Box$

\section{$L_1$-uniqueness}\label{S3}

In this section we prove Theorem~\ref{ttdep1.1}.
We adopt the Cauchy equation  approach of  Grigor'yan in his analysis of operators on  manifolds.
Grigor'yan's argument relies essentially on the geodesic completeness of the manifold but in the following proof this is 
replaced by  Markov uniqueness of $H$.
The latter property is equivalent, by Proposition~\ref{pcap2.1},  to $C_c^\infty(\Omega)$ being a core of $h_N$ and this
 suffices  for  the application of Grigor'yan's techniques.
 
 First for $\tau>0$ set $\Omega_\tau=\langle0,\tau\rangle\times \Omega$. 
 Denote a general point in $\Omega_\tau$ by
 $(t,x)$.
 So  $\partial_0$ denotes the partial derivative with respect to the first variable $t$.
A function $u\in L_\infty(\Omega_\tau)$ is defined to be a bounded weak solution of the 
 Cauchy equation corresponding to $K$ on $\Omega_\tau$ with initial value $\psi\in L_\infty(\Omega)$ if
 \begin{equation}
 (u, (-\partial_0+K)v)=0
 \label{ce1}
 \end{equation}
 for all $v\in C_c^\infty(\Omega_\tau)$ and 
 \begin{equation}
 \lim_{t\to0}\int_V dx\,|u(t,x)-\psi(x)|^2=0
 \label{ce2}
 \end{equation}
 for all relatively compact subsets $V$ of $\Omega$.
 Thus $u$ is a solution of the distributional equation $(-\partial_0+K)^*u=0$ on $\Omega_\tau$
 with initial condition $u(t,x)\to \psi(x)$ as $t\to 0$ in the $L_{2,\rm loc}(\Omega)$ sense.

The `time-dependent' criterion for  $L_1$-uniqueness of $K$  is 
formulated in terms of weak solutions of the  Cauchy equation with zero initial value.
 
 \begin{prop}\label{lcau3.1} If for some 
  $\tau>0$ the only bounded solution of  Cauchy equation $(\ref{ce1})$ on $\Omega_\tau$ with initial value $0$ in the $L_{2,\rm loc}$-sense $(\ref{ce2})$ is the zero solution
then $K$ is $L_1$-unique.
\end{prop}
\proof\
It follows from an extension of the Lumer--Phillips theorem (see \cite{Ebe}, Theorem~1.2 in Appendix~A of  Chapter~1) that $K$ is $L_1$-unique if and only if the $L_1$-closure of $K$ is the generator of a strongly continuous semigroup on $L_1(\Omega)$.
But this is the case if and only if  the range of $\lambda I+K$ is  $L_1$-dense for all large $\lambda>0$.

Assume that $K$ is not $L_1$-unique.
Thus for each large $\lambda$ there is a non-zero $\psi\in L_\infty(\Omega)$ such that $(\psi,(\lambda I+K)\varphi)=0$ for all  $\varphi\in C_c^\infty(\Omega)$.
Then  define $u_1$ on $\Omega_\tau$ by  $u_1(t,x)= e^{\lambda\,t}\,\psi(x)$ for all $t\in\langle0,\tau\rangle$ 
and all $x\in\Omega$.
It follows that $u_1$ is a  solution of the Cauchy equation (\ref{ce1}) on $\Omega_\tau$  with 
$\|u_1\|_\infty\leq e^{\lambda \tau}\|\psi\|_\infty$.
Moreover, $u_1$ has initial value $\psi$ in the $L_{2,\rm loc}$-sense (\ref{ce2}).

Next define $u_2$ on $\Omega_\tau$ by  $u_2(t,x)=(T_{t}^*\psi)(x)$ for all $t\in\langle0,\tau\rangle$ and $x\in\Omega$
where $T^*$ is the adjoint of the semigroup $T$ generated by the Friedrichs' extension $K_D$ of $K$.
The adjoint semigroup $T^*$ acts on $L_\infty(\Omega)$ and $\|T^*_s\|_{\infty\to\infty}=\|T_s\|_{1\to1}\leq e^{-(\omega_0-\omega_1)s}$ for all $s>0$ by the discussion of Subsection~\ref{S2.2}. 
Therefore $u_2$ is  also a  solution  of the Cauchy equation (\ref{ce1}) on $\Omega_\tau$ with 
$\|u_2\|_\infty\leq e^{\omega\tau}\|\psi\|_\infty$ where $\omega=(-\omega_0+\omega_1)\vee 0$.
But the adjoint semigroup $T^*$ on $L_\infty(\Omega)$ is  $L_{2, \rm loc}$-continuous by Proposition~\ref{pacc1}
and Remark~\ref{racc1}.
Thus $u_2$ has   initial value $\psi$ in the $L_{2,\rm loc}$-sense (\ref{ce2}).

Finally
\[
\sup_{x\in\Omega}|u_1(t,x)|\geq e^{(\lambda-\omega)t}\sup_{x\in\Omega}|u_2(t,x)|
\]
for all $t\in\langle0,\tau\rangle$.
Thus if $\lambda>\omega$ one must have $u_1\neq u_2$ and so $u_1-u_2$ is a non-zero bounded weak solution of the Cauchy equation (\ref{ce1}) with initial value zero in the $L_{2,\rm loc}$-sense (\ref{ce2}).

Therefore the proposition follows by negation.\hfill$\Box$
\bigskip

The key result in the proof of $L_1$-uniqueness, the analogue of Theorem~2 in \cite{Gri6}, Theorem~9.2 in  \cite{Gri5}
or Theorem~11.9 in \cite{Gri7},  can now be formulated as follows.
\begin{prop}\label{pnot1}
Assume $H$ is Markov unique and that the balls $B(r)$ are bounded for all $r>0$.
Let $u\in L_\infty(\Omega_\tau)$ be a bounded weak solution of the Cauchy equation $(\ref{ce1})$ with zero initial value in the  $L_{2,\rm loc}$-sense $(\ref{ce2})$.
Further assume 
\[
\int^\tau_0dt\int_{B(r)} dx\, |u(t,x)|^2\leq e^{\sigma(r)}
\]
for all large $r$ where $v$ is a positive increasing function on $\langle0,\infty\rangle$ such that 
\[
\int^\infty_R dr\,r\,\sigma(r)^{-1}=\infty
\]
for all large $R>0$.
Then $u=0$.
\end{prop}

This proposition in combination with Proposition~\ref{lcau3.1} immediately gives conditions for $L_1$-uniqueness
of $K$ or $H$. 

\begin{cor}\label{cnot1}
Assume the balls $B(r)$ are bounded for all $r>0$ and that the T\"acklind condition $(\ref{elcap1.31})$ is satisfied.
It follows that if $H$ is Markov unique then both $H$ and $K$ are $L_1$-unique.
\end{cor}
\proof\
Assume  $H$ is Markov unique.
If $u$ is a bounded weak solution of (\ref{ce1}) and $(\ref{ce2})$ then
\begin{equation}
\int^\tau_0dt\int_{B(r)} dx\,|u(t,x)|^2\leq \tau\,\|u\|_\infty^2\,|B(r)|
\label{ecau3.20}\;.
\end{equation}
It follows that the hypothesis of the proposition are fulfilled with $\sigma(r)=\log(\tau\,\|u\|_\infty^2\,|B(r)|)$.
Therefore $u=0$ by Proposition~\ref{pnot1} and $K$ is $L_1$-unique by Lemma~\ref{lcau3.1}.
But  setting the lower-order coefficients equal to zero one simultaneously deduces that $H$ is  $L_1$-unique.
\hfill$\Box$

\bigskip

The proof of Theorem~\ref{ttdep1.1}  is now reduced to proving Proposition~\ref{pnot1}.
Once this is established the theorem follows from Corollary~\ref{cnot1}.

\bigskip

\noindent{\bf Proof of Proposition~\ref{pnot1}} It suffices to prove that if  $r$ is large and $\delta\in\langle0,\tau]$
satisfies $\delta\leq  r^2/(16\,\sigma(r))$ then there is a $b>0$ such that 
\begin{equation}
\int_{B(r)}dx\,|u(\tau,x)|^2\leq \int_{B(2r)}dx\,|u(\tau-\delta,x)|^2+b\,r^{-2}\;.
\label{enot0}
\end{equation}
The rest of the proof then follows by direct repetition of  Grigoryan's argument \cite{Gri5} pages~186 and 187
or \cite{Gri7} pages~306 and 307.
In this  part of the proof, which we omit,  the $L_{2,\rm loc}$-initial condition is crucial.
Any weaker form of the initial condition is insufficient.
Now we concentrate on establishing~(\ref{enot0}).

\smallskip

Let $\rho_r(x)=\inf_{y\in B(r)}d(x\,;y)$ denote  the Riemannian distance from $x$ to the ball $B(r)$. 
Set $\xi_t=\nu\,\rho_r^2\,(t-s)^{-1}$ where $\nu, s>0$ are fixed with $t\neq s$.
The values of $s$ and $\nu$ will be chosen later. 
In particular the choice of $\nu$ depends on the lower-order coefficients.
It follows  that the partial derivative $\xi_t'$ with respect to $t$ is given by $\xi_t'=-\nu\,\rho_r^2\,(t-s)^{-2}$ and $\Gamma(\rho_r^2)=4\,\rho_r^2\,\Gamma(\rho_r)\leq4\,\rho_r^2$.
Therefore $\Gamma(\xi_t)\leq 4\,\nu^2\,\rho_r^2\,(t-s)^{-2}$ and 
\begin{equation}
\xi_t'+(4\,\nu)^{-1}\,\Gamma({\xi_t})\leq 0
\label{ecau3.3}
\;.
\end{equation}
(An auxiliary function of this type was introduced by Aronson, \cite{Aro} Section~3, in his derivation of Gaussian bounds on the heat kernel.)

First  we  consider the case that $u\in L_\infty(\Omega_\tau)$ is a weak solution of the Cauchy equation (\ref{ce1}) corresponding to $H$ and aim to deduce $L_1$-uniqueness of $H$. 
The argument for $K$ is very similar but the lower-order terms introduce additional computational complications.

In the notation of Subsection~\ref{S2.3} the Cauchy equation for $H$ states that $(u,\ch v)=0$ for all 
$v\in C_c^\infty(\Omega_\tau)$.
Therefore $u\in D(\ch^*)$. But $D(\ch^*)\subseteq V^{2,2}_{\rm loc}(\Omega_\tau)$ by Proposition~\ref{preg2.1}.
Thus the Cauchy equation can be explicitly written as
\begin{equation}
(\partial_0u,v)-\sum^d_{k,l=1}(\partial_kc_{kl}\partial_lu,v)=0
\label{ce3.1}
\end{equation}
for all $v\in C_c^\infty(\Omega_\tau)$.
But  (\ref{ce3.1}) extends to all $v\in L_2(\Omega_\tau)$ with compact support because  $u\in V^{2,2}_{\rm loc}(\Omega_\tau)$.
Now define $\psi_t$ by $\psi_t(x)=u(t,x)$ and let $\psi'_t$ denote its partial derivative with respect to $t$.
Then set $v$ equal to the restriction of $\eta^2e^{\xi_t}\psi_t$ to $\langle \tau-\delta,\tau\rangle\times\Omega$ with $\eta\in C_c^\infty(\Omega)$.
Thus $\supp v\subseteq[\tau-\delta, \tau]\times \supp\eta$ is compact.
It follows, after an integration by parts in the $x$-variables,  that 
\begin{eqnarray}
\int^\tau_{\tau-\delta}dt\,(\psi_t', \eta^2e^{2\xi_t}\psi_t)&=&-\sum^d_{k,l=1}\int^\tau_{\tau-\delta}dt\,(\partial_l\psi_t, c_{kl}\,\partial_k(\eta^2\,e^{2\xi_t}\psi_t))\nonumber\\[-8pt]
&=&\int^\tau_{\tau-\delta}dt\,(\psi_t, \Gamma(\eta e^{\xi_t})\psi_t)-\sum^d_{k,l=1}\int^\tau_{\tau-\delta}dt\,(\partial_l(\eta e^{\xi_t}\psi_t), c_{kl}\partial_k(\eta e^{\xi_t}\psi_t))
\nonumber\\[-1pt]
&=&\int^\tau_{\tau-\delta}dt\,(\psi_t, \Gamma(\eta e^{\xi_t})\psi_t)-\int^\tau_{\tau-\delta}dt\,h_D(\eta e^{\xi_t}\psi_t)
\;.
\label{enot20}
\end{eqnarray}
Since $\eta\in C_c^\infty(\Omega)$ there are no boundary terms.
But one also has \begin{eqnarray}
(\psi_t, \Gamma(\eta e^{\xi_t})\psi_t)&\leq & 2\,(e^{\xi_t}\psi_t,\Gamma(\eta)e^{\xi_t}\psi_t)+2\,(\eta\psi_t,\Gamma(e^{\xi_t})\eta\psi_t)\nonumber\\[5pt]
&= &2\,(e^{\xi_t}\psi_t,\Gamma(\eta)e^{\xi_t}\psi_t)+2\,(\eta e^{\xi_t}\psi_t,\Gamma(\xi_t)\eta e^{\xi_t}\psi_t)
\label{enot201}
\end{eqnarray}
for all $\eta\in C_c^\infty(\Omega)$.
Combination of (\ref{enot20})
 and (\ref{enot201}) immediately leads to the  inequality
\begin{eqnarray*}
2^{-1}\int^\tau_{\tau-\delta}dt\int_\Omega \eta^2e^{2\xi_t}
(\psi_t^2)'
&\leq &2\int^\tau_{\tau-\delta}dt\,(e^{\xi_t}\psi_t,\Gamma(\eta)e^{\xi_t}\psi_t)+2\int^\tau_{\tau-\delta}dt\,
(\eta e^{\xi_t}\psi_t,\Gamma(\xi_t)\eta e^{\xi_t}\psi_t)
\end{eqnarray*}
for all $\eta\in C_c^\infty(\Omega)$.
Then integrating by parts in the $t$-variable and rearranging gives
\begin{eqnarray}
\Big[\|\eta \,e^{\xi_t}\psi_t\|_2^2\Big]^\tau_{\tau-\delta}
&\leq&2\int^\tau_{\tau-\delta}dt\,(e^{\xi_t}\psi_t,\Gamma(\eta)e^{\xi_t}\psi_t)+
\int^\tau_{\tau-\delta}dt\,(\eta e^{\xi_t}\psi_t,(\xi_t'+2\,\Gamma(\xi_t))\eta e^{\xi_t}\psi_t)\nonumber\\[5pt]
&\leq &2\int^\tau_{\tau-\delta}dt\,(e^{\xi_t}\psi_t,\Gamma(\eta)e^{\xi_t}\psi_t)\label{enot21}
\end{eqnarray}
for all $\eta\in C_c^\infty(\Omega)$
where the last step uses (\ref{ecau3.3}) with $\nu$ chosen equal to $8^{-1}$.
Next we use the Markov uniqueness of $H$ to extend (\ref{enot21}) to a larger class of $\eta$.

First  choose $s=\tau+\delta$ in the definition of $\xi_t$ so  with the previous choice of $\nu=8^{-1}$ one has $\xi_t=-8^{-1}\rho_r^2(\tau+\delta-t)^{-1}\leq0$ for all $t\in\langle0,\tau]$.
Therefore 
\[
\|\eta \,e^{\xi_t}\psi_t\|_2\leq \|\eta \,\psi_t\|_2\leq \|\psi_t\|_\infty \|\eta\|_2\leq \|u\|_\infty\|\eta\|_2
\]
and 
\[
0\leq (e^{\xi_t}\psi_t,\Gamma(\eta)e^{\xi_t}\psi_t)\leq  (\psi_t,\Gamma(\eta)\psi_t)\leq \|u\|_\infty\,h_D(\eta)
\]
for all $\eta\in C_c^\infty(\Omega)$ and all $t\in\langle0,\tau]$.
Since $H$ is Markov unique $h_D=h_N$ and $C_c^\infty(\Omega)$ is a core of $h_N$ by Proposition~\ref{pcap2.1}.
It then follows 
by continuity that (\ref{enot21}) extends  to all $\eta\in D(h_N)$.
Thus one concludes that 
\begin{equation}
\|\eta \,e^{\xi_\tau}\psi_\tau\|_2^2\leq \|\eta \,e^{\xi_{\tau-\delta}}\psi_{\tau-\delta}\|_2^2
+2\int^\tau_{\tau-\delta}dt\,(e^{\xi_t}\psi_t,\Gamma(\eta)e^{\xi_t}\psi_t)\label{enot22}
\end{equation}
for all $\eta\in D(h_N)$.

Next  let $\theta\in C_c^\infty(\Ri)$ satisfy  $0\leq \theta\leq1$, $\theta(s)=1$ if $s\in[0,3/2]$,  $\theta(s)=0$ if $s\geq 2$ and $|\theta'|\leq 3$. 
Then  set $\theta_r=\theta\circ(r^{-1}\rho)$.
It follows that  $\theta_r\in D(h_N)\cap L_\infty(\Omega)$. 
Moreover, $\theta_r=1$ if $\rho\leq 3r/2$ and $\theta_r=0$ if $\rho\geq 2r$.
Thus $\supp\theta_r\subseteq B(2r)$ which is a bounded subset of $\Omega$ by assumption.
But  $\Gamma(\rho)\leq 1$.
So  one also has  $\|\Gamma(\theta_r)\|_\infty\leq 9\,r^{-2}$.
Hence replacing $\eta$ in  (\ref{enot22}) by $\theta_r$ one has
\begin{equation}
\int_{B(r)}|e^{\xi_\tau}\psi_\tau|^2\leq \int_{B(2r)}|e^{\xi_{\tau-\delta}}\psi_{\tau-\delta}|^2
+18(a/r)^2\int^\tau_{\tau-\delta}dt\int_{B(2r)\backslash B(3r/2)}dx\,|(e^{\xi_t}\psi_t)(x)|^2
\label{enot3}
\end{equation}
But if $x\in B(r)$ then $\xi_\tau=0$.
Moreover, $\xi_{\tau-\delta}\leq 0$.
Further if $x\in B(2r)\backslash B(3r/2)$ then $\rho_r(x)\geq r/2$ and so $\xi_t(x)\leq -r^2/(16\delta)$ for $t\in \langle \tau-\delta,\tau\rangle$.
Then it follows from (\ref{enot3})  and the hypothesis of the proposition that 
\begin{eqnarray}
\int_{B(r)}|\psi_\tau|^2&\leq& \int_{B(2r)}|\psi_{\tau-\delta}|^2
+18(a/r)^2\int^\tau_{\tau-\delta}dt\int_{B(2r)}dx\,|\psi_t|^2e^{-r^2/(16\delta)}\nonumber\\[5pt]
&\leq& \int_{B(2r)}|\psi_{\tau-\delta}|^2
+18(a/r)^2e^{-(r^2/(16\delta))+\sigma(2r)}\;.
\label{enot4}
\end{eqnarray}
Finally choosing $\delta\leq r^2/(16\sigma(2r))$ one has
\[
\int_{B(r)}|\psi_\tau|^2\leq \int_{B(2r)}|\psi_{\tau-\delta}|^2+18(a/r)^2
\;.
\]
Thus  we have established (\ref{enot0}) and  the proposition follows for a solution of the Cauchy equation corresponding to $H$.
Thus $H$ is $L_1$-unique.

\smallskip

In order to conclude that $K$ is $L_1$-unique it remains to prove Proposition~\ref{pnot1}
for a solution of the Cauchy equation (\ref{ce1})  corresponding to $K=H+L+M$.
In particular  we have to consider the estimation of the lower-order terms.
But now with the notation of Subsection~\ref{S2.3} the Cauchy equation  states that
\[
(u,\ck_0v)+(u,Mv)=0
\]
for all $v\in C_c^\infty(\Omega_\tau)$.
Let $V_\tau=\langle0,\tau\rangle\times V$ where $V$ is a relatively compact subset of $\Omega$.
It follows that
\[
|(u,\ck_0v)|\leq \|u\|_\infty\|Mv\|_1\leq \tau^{1/2}\|u\|_\infty\|c_0\|_{L_2(V)}\|v\|_2
\]
for all $v\in C_c^\infty(V_\tau)$ because of the assumption that $c_0\in L_{2,\rm loc}(\Omega)$.
Hence  $u$ is in the domain of the adjoint of $\ck_0|_{C_c^\infty(V_\tau)}$.
Then one deduces from Proposition~\ref{preg2.2} that $u\in V^{2,2}_{\rm loc}(\Omega_\tau)$.
Therefore one can argue as before.
First the Cauchy equation (\ref{ce3.1}) is replaced by
\begin{equation}
(\partial_0u,v)-\sum^d_{k,l=1}(\partial_kc_{kl}\partial_lu,v)+(u,Lv)+(u,Mv)=0
\label{ce3.11}
\end{equation}
for all $v\in C_c^\infty(\Omega_\tau)$.
Then (\ref{enot20}) is replaced by

\begin{eqnarray}
\int^\tau_{\tau-\delta}dt\,(\psi_t', \eta^2e^{2\xi_t}\psi_t)&=&\int^\tau_{\tau-\delta}dt\,(\psi_t, \Gamma(\eta e^{\xi_t})\psi_t)-\int^\tau_{\tau-\delta}dt\,h_D(\eta e^{\xi_t}\psi_t)\nonumber\\[5pt]
&&\hspace{0.5cm}{}-\int^\tau_{\tau-\delta}dt\,(\psi_t, L\eta^2e^{2\xi_t}\psi_t)
-\int^\tau_{\tau-\delta}dt\,(e^{\xi_t}\eta\psi_t,Me^{\xi_t}\eta\psi_t)
\;.
\label{enot40}
\end{eqnarray}
The first term on the right hand side is again estimated by (\ref{enot201}) and it remains to estimate the terms originating with the lower-order terms $L$ and $M$.
But 
\begin{eqnarray*}
(\psi_t, L\eta^2e^{2\xi_t}\psi_t)
=(\eta e^{\xi_t}\psi_t, L\eta e^{\xi_t}\psi_t)+(\psi_t,[L,\eta e^{\xi_t}]\eta e^{\xi_t}\psi_t)
\end{eqnarray*}
Further
\begin{eqnarray*}
(\psi_t,[L,\eta e^{\xi_t}]\eta e^{\xi_t}\psi_t)
=(e^{\xi_t}\psi_t,\eta L(\eta)e^{\xi_t}\psi_t)+(e^{\xi_t}\eta\psi_t,L(\xi_t)\eta e^{\xi_t}\psi_t)
\end{eqnarray*}
It follows, however, from the third condition in (\ref{elcap1.12}) that
\[
\|L(\eta)\varphi\|_2^2\leq \kappa^{-1}(\varphi, \Gamma(\eta)\varphi)
\]
for all $\varphi\in L_2(\Omega)$.
Therefore
\begin{equation}
\left.
\begin{array}{rl}
|(\eta e^{\xi_t}\psi_t, L\eta e^{\xi_t}\psi_t)|\hspace{-2mm}&\leq h_D(\eta e^{\xi_t}\psi_t)+(4\kappa)^{-1}(\eta e^{\xi_t}\psi_t,\eta e^{\xi_t}\psi_t)\;,
\\[10pt]
|(e^{\xi_t}\psi_t,\eta L(\eta)e^{\xi_t}\psi_t)|\hspace{-2mm}&\leq 2^{-1}\,(\eta e^{\xi_t}\psi_t,\eta e^{\xi_t}\psi_t)
+(2\kappa)^{-1}(e^{\xi_t}\psi_t,\Gamma(\eta) e^{\xi_t}\psi_t)
\;,\\[10pt]
|(e^{\xi_t}\eta\psi_t,L(\xi_t)\eta e^{\xi_t}\psi_t)|\hspace{-2mm}&\leq 2^{-1}\,(\eta e^{\xi_t}\psi_t,\eta e^{\xi_t}\psi_t)+
(2\kappa)^{-1}(\eta e^{\xi_t}\psi_t,\Gamma(\xi_t) \eta e^{\xi_t}\psi_t)
\;.
\end{array}
\right\}\label{enot401}
\end{equation}
Combining   estimates (\ref{enot201}),  (\ref{enot40}) and (\ref{enot401}) one deduces that 
\begin{eqnarray*}
(\psi_t', \eta^2e^{2\xi_t}\psi_t)
&\leq &2\gamma\,(e^{\xi_t}\psi_t,\Gamma(\eta)e^{\xi_t}\psi_t)
+2\gamma\,(\eta e^{\xi_t}\psi_t,\Gamma(\xi_t)\eta e^{\xi_t}\psi_t)\\[5pt]
&&\hspace{1.5cm}{}-(\eta e^{\xi_t}\psi_t,(c_0-1-\gamma)\eta e^{\xi_t}\psi_t)
\end{eqnarray*}
with $\gamma=(1+(4\kappa)^{-1})$.
Now $L_1$-uniqueness of $K$ is equivalent to $L_1$-uniqueness of $K+\omega I$ for any $\omega\in \Ri$.
Therefore, replacing $c_0$ by $c_0+\omega$  one may assume $\omega_0\geq 1+ \gamma$.
Hence $c_0-1-\gamma\geq0$ and one concludes that
\begin{eqnarray*}
(\psi_t', \eta^2e^{2\xi_t}\psi_t)
\leq 2\gamma\,(e^{\xi_t}\psi_t,\Gamma(\eta)e^{\xi_t}\psi_t)
+2\gamma\,(\eta e^{\xi_t}\psi_t,\Gamma(\xi_t)\eta e^{\xi_t}\psi_t)
\;.
\end{eqnarray*}
Integrating by parts and rearranging gives
\begin{eqnarray}
\Big[\|\eta \,e^{\xi_t}\psi_t\|_2^2\Big]^\tau_{\tau-\delta}
&\leq&2\gamma\int^\tau_{\tau-\delta}dt\,(e^{\xi_t}\psi_t,\Gamma(\eta)e^{\xi_t}\psi_t)\nonumber\\[5pt]
&&\hspace{2cm}{}+\int^\tau_{\tau-\delta}dt\,(\eta e^{\xi_t}\psi_t,(\xi_t'+2\gamma\,\Gamma(\xi_t))\eta e^{\xi_t}\psi_t)\end{eqnarray}
for all $\eta\in C_c^\infty(\Omega)$.
Then setting  $\nu=(8\gamma)^{-1}$ in the definition of $\xi_t$ one has
$\xi_t'+2\gamma\,\Gamma(\xi_t)\leq0$
by (\ref{ecau3.3}). 
Therefore one concludes that 
\begin{equation}
\|\eta \,e^{\xi_\tau}\psi_\tau\|_2^2\leq \|\eta \,e^{\xi_{\tau-\delta}}\psi_{\tau-\delta}\|_2^2
+2\gamma\int^\tau_{\tau-\delta}dt\,(e^{\xi_t}\psi_t,\Gamma(\eta)e^{\xi_t}\psi_t)\label{enot42}
\end{equation}
for all $\eta\in C_c^\infty(\Omega)$ in direct analogy with (\ref{enot21}).
In particular this estimate is valid with $s=\tau+\delta$ in the definition of $\xi_t$.
Since $H$ is Markov unique (\ref{enot42})  extends to all $\eta\in D(h_N)$
by repetition of the previous reasoning.
The rest of the proof   is exactly the same as the earlier proof  for $H$.
Using (\ref{enot42}) in place of (\ref{enot22}) one establishes Proposition~\ref{pnot1} for $K$ and thereby concludes that $K$
is $L_1$-unique. \hfill$\Box$

\bigskip

The foregoing `time-dependent' argument to deduce $L_1$-uniqueness from Markov uniqueness appears to be quite different
to the `time-independent' arguments of \cite{RSi5} and \cite{RSi4} for the symmetric operator $H$.
The two methods are, however,  related.
The time-independent proof uses Davies--Gaffney off-diagonal Gaussian bounds \cite{Gaf} \cite{Dav12}
and  one derivation  of the latter bounds is by a variation of the foregoing time-dependent argument. 
(See  \cite{Gri5} Chapter~12.)
The time-dependent argument is based on the T\"acklind condition (\ref{elcap1.31}) on $|B(r)|$ but
 the time-independent method for $H$ requires the stronger condition $|B(r)|\leq a\,e^{b\,r^2}$ for some $a,b>0$ and all $r>0$.
The latter restriction is essential because the  argument uses the Davies--Gaffney off-diagonal bounds.

One may extend Theorem~\ref{ttdep1.1} to operators $K$ for which the coefficients $c_k$ and $c_0$ are complex-valued.
But then the assumptions (\ref{elcap1.12}) have to be appropriately modified, e.g.\ it is necessary that 
$\RRe c_0$ is lower semi-bounded
and $\RRe\divv c$ is upper semi-bounded,
Moreover, the third condition in (\ref{elcap1.12}) has to be replaced by 
$C(x)\geq \kappa\, \big(\,\overline c^T\otimes c+c^T\otimes\overline c\,\big)(x) $
for almost all $x\in\Omega$.
The proof is essentially the same but the spaces involved are complex.

\section{Markov uniqueness}
\label{S5}

The   basic ingredients in the foregoing analysis of $L_1$-uniqueness were the growth restrictions on the Riemannian geometry and the  Markov uniqueness of $H$.
In this section we  consider the characterization of  the latter property by capacity conditions. 
The first result of this nature is due to  Maz'ya (see \cite{Maz} Section~2.7) for the case $\Omega=\Ri^d$.
Maz'ya demonstrated that the identity  $h_D=h_N$  is equivalent to a family of conditions
on sets of finite capacity.
More recently it was established in \cite{RSi5} and \cite{RSi4} that Markov uniqueness is equivalent to the capacity of the boundary of $\Omega$ being zero.
Our aim is to  establish that both these capacity criteria are valid for $H$ and for general open $\Omega$
whenever the Riemannian balls $B(r)$ are bounded for all $r>0$.
But this requires in part a slightly stronger assumption on the properties of the coefficients $c_{kl}$.

First we define a  subset $A$ of $\overline\Omega$  to have finite capacity, relative to $H$, if there is an $\eta\in D(h_N)$ such that $\eta=1$ on $A$.
Each relatively compact subset of $\Omega$ has finite capacity by Urysohn's lemma.
Moreover, each set of finite capacity $A$ must have finite volume, i.e.\ $|A|<\infty$, but 
one can have  unbounded sets with finite capacity (see, \cite{Maz}, Section~2.7).

We begin by establishing that there are an abundance of sets of finite capacity.

\begin{prop}\label{pcap2.11} The subspace $(D(h_N)\cap L_\infty(\Omega))_{\rm cap}$ of bounded functions in $D(h_N)$
whose supports have finite capacity is a core of $h_N$.
\end{prop}
\proof\ 
It suffices to prove that each $\varphi\in D(h_N)\cap L_\infty(\Omega)$ can be approximated in the $D(h_N)$-graph norm by a sequence 
 $\varphi_n\in (D(h_N)\cap L_\infty(\Omega))_{\rm cap}$.
 Clearly one may assume that $\varphi\geq0$.
But  if $\lambda>0$  the set $A_\lambda=\{x\in \Omega:\varphi(x)> \lambda\}$  has  finite capacity.
 This is a consequence  of the Dirichlet form structure by the following argument of Maz'ya.
 Define $\varphi_\lambda$ by $\varphi_\lambda(x)=\lambda^{-1}(\varphi(x)\wedge \lambda)$.
 Then $\varphi_\lambda\in D(h_N)$,  $0\leq \varphi_\lambda\leq 1$,  $\varphi_\lambda=1$ on $A_\lambda$ and 
 $h_N(\varphi_\lambda)\leq \lambda^{-2}h_N(\varphi)$  where the latter bounds follows from  the Dirichlet property of $h_N$.
 Therefore $A_\lambda$ has finite capacity.
 Now consider the sequence  $\varphi_m=\varphi-\varphi\wedge m^{-1}\in D(h_N)\cap L_\infty(\Omega)$.
 Since  $\supp\varphi_m= A_{m^{-1}}$ it follows that  $\varphi_m\in (D(h_N)\cap L_\infty(\Omega))_{\rm cap}$.
But the   $\varphi_m$  converge in the $D(h_N)$-graph norm to $\varphi$
as $m\to\infty$  by  \cite{FOT}, Theorem~1.4.2(iv).
\hfill$\Box$

\bigskip

Secondly,  to formulate suitable versions of Mazya's approximation criterion for Markov uniqueness we introduce the condition $\cc_A$ for each subset $A$ of 
$\overline\Omega$ by
 \begin{equation}
\hspace{-0cm}\mbox{$\cc_A$:}\hspace{0.8cm}
 \left\{ \begin{array}{ll}
\hspace{0mm}\mbox{ there exist }
\eta_1,\eta_2, \ldots \in D(h_D)
\mbox{ such that }&{}\\[8pt]
\hspace{6mm}\lim_{n\to\infty}\|\one_A\,(\one_\Omega-\eta_n )\, \varphi\|_2 = 0
\mbox{ for each } \varphi \in L_2(\Omega)&{}\\[8pt]
\hspace{0mm}\mbox{ and } \lim_{n\to\infty} \|\one_A\,\Gamma(\eta_n)\|_1 = 0.&{}
  \end{array}\right.\label{gc}
\end{equation}
Although  the approximating sequence in this condition is  formed by functions $\eta_n\in D(h_D)$ one can, equivalently,
choose  $\eta_n\in C_c^\infty(\Omega)$. 
This follows because   $C_c^\infty(\Omega)$ is a core of $h_D$.
Explicitly, for each $\eta_n\in D(h_D)$ there is a $\chi_n\in C_c^\infty(\Omega)$ such that $\|\eta_n-\chi_n\|_{D(h_D)}\leq n^{-1}$.
Therefore
\[
\|\one_A\,(\one_\Omega-\chi_n )\, \varphi\|_2\leq \|\one_A\,(\one_\Omega-\eta_n )\, \varphi\|_2+n^{-1}\|\varphi\|_\infty
\]
for all $\varphi\in L_2(\Omega)\cap L_\infty(\Omega)$ and
\[
\|\one_A\,\Gamma(\chi_n)\|_1\leq \|\one_A\,\Gamma(\eta_n)\|_1+h_D(\chi_n-\eta_n)
\leq \|\one_A\,\Gamma(\eta_n)\|_1+n^{-1}
\;.
\]
Hence the $\cc_A$-convergence criteria for the $\chi_n$ are inherited from the $\eta_n$.
Alternatively, one may assume, without loss of generality, that the $\eta_n$ satisfy $0\leq \eta_n\leq 1$.
This follows because $\zeta_n=(0\vee \eta_n)\wedge 1\in D(h_N)$, 
\[
\|\one_A\,(\one_\Omega-\zeta_n )\, \varphi\|_2\leq \|\one_A\,(\one_\Omega-\eta_n )\, \varphi\|_2
\]
for all $\varphi\in L_2(\Omega)$ and $\|\one_A\,\Gamma(\zeta_n)\|_1\leq \|\one_A\,\Gamma(\eta_n)\|_1$
(see \cite{BH}, Proposition~4.1.4).
Therefore the $\zeta_n$ inherit the $\cc_A$-convergence properties of the $\eta_n$.

The next proposition is a local version of Maz'ya's result \cite{Maz}, Theorem~1 in Section~2.7 (see also \cite{FOT}, Theorem 3.2.2).
Note that it is independent of any constraints on the Riemannian geometry.
\begin{prop}\label{pmaz}
The following conditions are equivalent:
\begin{tabel}
\item\label{pmaz1}
$H$ is Markov unique,
\item\label{pmaz2}
$\cc_A$ is satisfied for each subset $A$ of $\overline\Omega$ with finite capacity.
\end{tabel}
\end{prop}
\proof\ \ref{pmaz1}$\Rightarrow$\ref{pmaz2}$\;$
If  $A\subseteq\overline\Omega$ is  a set of finite capacity there exists an  $\eta\in D(h_N)$ with  $\eta=1$ on~$A$.
But $h_N=h_D$, by  Markov uniqueness.
Therefore $\eta\in D(h_D)$.
Then the constant sequence $\eta_n=\eta$ satisfies $\cc_A$.

\smallskip

\noindent\ref{pmaz2}$\Rightarrow$\ref{pmaz1}$\;$
It suffices to prove that each $\varphi\in D(h_N)$ can be approximated in the $D(h_N)$-graph norm by a sequence 
 $\varphi_n\in D(h_D)\cap L_\infty(\Omega)$.
 But $ (D(h_N)\cap L_\infty(\Omega))_{\rm cap}$ is a core of $h_N$, by Proposition~\ref{pcap2.11}. 
 Therefore one may assume that $ \varphi\in(D(h_N)\cap L_\infty(\Omega))_{\rm cap}$.
 Set $A=\supp\varphi$ and let $\eta_n\in D(h_D)\cap L_\infty(\Omega)$ be the  corresponding  $\cc_A$-sequence.
Then let $\varphi_n=\eta_n\varphi$. 
 It follows from Proposition~\ref{pcap2.10} that $\varphi_n\in D(h_D)\cap L_\infty(\Omega)$. 
But 
 \[
\lim_{n\to\infty} \|\varphi-\varphi_n\|_2=\lim_{n\to\infty}\|\one_A(\one_\Omega-\eta_n)\,\varphi\|_2=0
\;.
 \]
In addition $\nabla(\varphi_n-\varphi)=(\nabla\eta_n)\,\varphi+(1-\eta_n)\,(\nabla\varphi)$.
Therefore
\[
\Gamma(\varphi_n-\varphi)\leq 2\,\Gamma(\eta_n)\,\varphi^2+2\,(1-\eta_n)^2\,\Gamma(\varphi)
\;.
\]
Then since $\supp\varphi_n\subseteq\supp\varphi$ it follows that
\[
h_N(\varphi-\varphi_n)=\|\one_A\Gamma(\varphi_n-\varphi)\|_1\leq 
2\,\|\one_A\Gamma(\eta_n)\|_1\,\|\varphi\|_\infty^2+2\,\|\one_A(\one_\Omega-\eta_n)\chi\|_2^2
\]
where $\chi=\Gamma(\varphi)^{1/2}\in L_2(\Omega)$.
Therefore $h_N(\varphi-\varphi_n)\to0$ as $n\to\infty$.
This establishes
 that $D(h_D)\cap L_\infty(\Omega)$ is a core of $h_N$.
 Hence $h_D=h_N$ and $H$ is Markov unique. \hfill$\Box$

 \bigskip

 Next we discuss improvements to the foregoing results with two additional assumptions.
 First we assume the Riemannian balls $B(r)$ are bounded for all $r>0$.
 This  immediately gives an improved version of Proposition~\ref{pcap2.11}.

\begin{prop}\label{pcap2.2} 
Assume $B(r)$ is bounded for all $r>0$.
Then $(D(h_N)\cap L_\infty(\Omega))_c$,  the subspace of bounded functions in $D(h_N)$
with compact support, is a core of $h_N$.
\end{prop}
\proof\
The proof is essentially identical to the proof of Lemma~2.3 in \cite{OuR}.
First,   the $B(r)$ are bounded if and only if $\rho(x)\to\infty$ as $x\to\infty$ where $\rho$ is again the Riemannian distance from the origin.
Secondly, let $\tau\in C_c^\infty(\Ri)$ satisfy  $0\leq \tau\leq1$, $\tau(s)=1$ if $s\in[0,1]$,  $\tau(x)=0$ if $s\geq 2$ and $|\tau'|\leq 2$. 
Then  set $\tau_n=\tau\circ(n^{-1}\rho)$.
It follows that  $\tau_n$ has compact support. 
Moreover,
 $\tau_n(x)\to1$ as $n\to\infty$ for all $x\in\Omega$.
But  $\Gamma(\rho)\leq 1$.
So  one also has  $\|\Gamma(\tau_n)\|_\infty\leq 4\,n^{-2}$.
Thirdly,  if $\varphi\in D(h_N)\cap L_\infty(\Omega)$ and  $\varphi_n=\tau_n\,\varphi$ then $\varphi_n\in (D(h_N)\cap  L_\infty(\Omega))_c$ by Proposition~\ref{pcap2.10}.
But 
\[
\|\varphi_n-\varphi\|_{D(h_N)}^2\leq 
2\int_{\Omega}\Gamma(\tau_n)\,\varphi^2
+2\int_\Omega(\one_\Omega-\tau_n)^2\,\Gamma(\varphi)
+ \int_{\Omega}(\one_\Omega-\tau_n)^2\,\varphi^2
\]
and all three terms on the right converge to zero as  $n\to\infty$ by the dominated convergence theorem.
Therefore $\varphi\in D(h_N)\cap L_\infty(\Omega)$ is the limit 
of the  $\varphi_n\in (D(h_N)\cap  L_\infty(\Omega))_c$ with respect to the $D(h_N)$-graph norm.
 Since $D(h_N)\cap L_\infty(\Omega)$ is a core of $h_N$ it follows that $(D(h_N)\cap  L_\infty(\Omega))_c$ is also a core.
\hfill$\Box$

  \bigskip

Secondly we assume that the  $c_{kl}\in W^{1,\infty}_{\rm loc}(\overline\Omega)$, the space of restrictions to $\Omega$
of functions in $W^{1,\infty}_{\rm loc}(\Ri^d)$.
This  ensures that the coefficients extend by continuity to functions which are uniformly locally bounded
on $\overline\Omega$.
Therefore each bounded subset $A$ of $\overline\Omega$ has finite capacity.
This is again a consequence of Urysohn's lemma.
Now one can establish  an improved version of Proposition~\ref{pcap2.11}.

\begin{thm}\label{ntcap2.11}
Assume $c_{kl}\in W^{1,\infty}_{\rm loc}(\overline\Omega)$.
Consider the following conditions:
\begin{tabel}
\item\label{ntcap2.11-3}
$H$ is Markov unique,
\item\label{ntcap2.11-40}
$\cc_A$ is satisfied for each bounded 
subset $A$ of $\overline \Omega$.
\end{tabel}
Then 
{\rm \ref{ntcap2.11-3}$\Rightarrow$\ref{ntcap2.11-40}}.
Moreover, if $B(r)$ is bounded for all $r>0$  then
 {\rm \ref{ntcap2.11-40}$\Rightarrow$\ref{ntcap2.11-3}} and the conditions are equivalent.
\end{thm}
\proof\ \noindent \ref{ntcap2.11-3}$\Rightarrow$\ref{ntcap2.11-40}$\;$%
If $A\subseteq\overline\Omega$  is bounded then there is an $\eta\in C_c^\infty(\Ri^d)$ with  $\eta=1$ on~$A$.
Since the $c_{kl}\in W^{1,\infty}_{\rm loc}(\overline\Omega)$
it follows that $\eta$, or more precisely the restriction of $\eta$ to $\Omega$, is in $D(h_N)$.
But $h_N=h_D$, by  Markov uniqueness. 
Therefore $\eta\in D(h_D)$ and  the constant sequence $\eta_n=\eta$ satisfies Condition $\cc_A$.

\smallskip

This establishes the first statement in Theorem~\ref{ntcap2.11}.
Next we   assume that the balls $B(r)$ are bounded and consider the converse reasoning.

\smallskip

\noindent \ref{ntcap2.11-40}$\Rightarrow$\ref{ntcap2.11-3}$\;$It is necessary to prove that $D(h_N)=D(h_D)$.
 But  since the Riemannian balls $B(r)$ are bounded $(D(h_N)\cap L_\infty(\Omega))_c$ is a core of $h_N$ by Proposition~\ref{pcap2.2}. 
Therefore it suffices to prove that each $\varphi\in (D(h_N)\cap L_\infty(\Omega))_c$
 can be approximated in the $D(h_N)$-graph norm by a sequence 
 $\varphi_n\in D(h_D)\cap L_\infty(\Omega)$.

 Let $A=\supp\varphi$. 
  If $\eta_n\in D(h_D)$ is the $\cc_A$-sequence corresponding to the bounded set $A$ define $\varphi_n$
 by $\varphi_n=\eta_n\,\varphi$.
 Since we may assume $\eta_n\in D(h_D)\cap L_\infty(\Omega)$ it follows that $\varphi_n\in D(h_D)\cap L_\infty(\Omega)$ 
 by Proposition~\ref{pcap2.10}.
 But then the argument used to prove \ref{pmaz2}$\Rightarrow$\ref{pmaz1} in Proposition~\ref{pmaz}
 establishes that $\varphi_n$ converges to $\varphi$ in the $D(h_N)$-graph norm.
Therefore $D(h_D)\cap L_\infty(\Omega)$ is a core of $h_N$.
 Hence $h_D=h_N$ and $H$ is Markov unique.
 \hfill$\Box$

 \bigskip
 
 The assumption that the balls $B(r)$ are bounded is essential for the implication \ref{ntcap2.11-40}$\Rightarrow$\ref{ntcap2.11-3} in Theorem~\ref{ntcap2.11}.
 Maz'ya  has constructed an example for  $\Omega=\Ri^d$ (see, \cite{Maz}, Theorem~3 in Section~2.7) 
in which the  coefficients  grow rapidly in a  set with an infinitely extended cusp.
 The growth is such that  the Riemannian distance to infinity along the axis of the cusp is finite 
and consequently the balls $B(r)$ are not bounded for all sufficiently large $r$.
 In this  example   Condition~\ref{ntcap2.11-40} 
 is satisfied
but $h_D\neq h_N$, i.e.\  Condition~\ref{ntcap2.11-3} is false. 

 \smallskip

 Condition~\ref{ntcap2.11-40} of Theorem~\ref{ntcap2.11} is related to the boundary capacity condition established in  Theorem~1.2 in \cite{RSi5} as a characterization of Markov uniquenes.
 We  conclude this section with a brief discussion of the relationship.
The capacity of  a general subset $A$ of $\overline\Omega$ is defined by
\begin{eqnarray*}
\capp(A)=\inf\Big\{\;\|\psi\|_{D(h_N)}^2&&\;: \;\psi\in D(h_N), \;0\leq\psi\leq 1  \mbox{ and  there exists   an open set  }\nonumber\\[-5pt]
&& U\subset \Ri^d
\mbox{ such that } U\supseteq A
\mbox{ and } \psi=1 \mbox{ on } U\cap\Omega\;\Big\}
\end{eqnarray*}
with the convention that $\capp(A)=\infty$ if the infimum is over the empty set.
(This definition is analogous  to the  canonical definition of the capacity associated with a Dirichlet form
\cite{BH}  \cite{FOT}
and  if $\Omega=\Ri^d$ the two definitions coincide.) 
A slight extension of the arguments of \cite{RSi5} then gives the following characterization of Markov uniqueness.

\begin{prop}\label{pcap}
Assume $c_{kl}\in W^{1,\infty}_{\rm loc}(\overline\Omega)$. Consider the following conditions:
\begin{tabel}
\item\label{pcap1}
$H$ is Markov unique,
\item\label{pcap2}
$\capp(\partial\Omega)=0$.
\end{tabel}
Then {\rm \ref{pcap1}$\Rightarrow$\ref{pcap2}}.
Moreover, if  the balls $B(r)$ are bounded for all $r>0$  then
{\rm \ref{pcap2}$\Rightarrow$\ref{pcap1}} and the conditions are equivalent.
\end{prop}
\proof\ \ref{pcap1}$\Rightarrow$\ref{pcap2}$\;$  It follows from the general  properties of the capacity that if $A=\bigcup_{k=1}^\infty A_k$ then $\capp(A)\leq \sum_{k=1}^\infty\capp(A_k)$.
Therefore it  suffices to prove that $\capp(B)=0$ for each bounded $B\subseteq\partial\Omega$.
But  $\capp(B)<\infty$, because $c_{kl}\in W^{1,\infty}_{\rm loc}(\overline\Omega)$.
Therefore there is an open subset  $U$ of $\Ri^d$ containing $B$ and a $\psi\in D(h_N)$ with $\psi=1$ on $U\cap\Omega$. 
Then by  Markov uniqueness 
one can  find  $\psi_n\in (D(h_D)\cap L_\infty(\Omega))_c$ such that  $\|\psi-\psi_n\|_{D(h_N)}\to0$ as $n\to\infty$.
Therefore there are open subsets $U_n$ of $\Ri^d$ 
 containing $B$  with $\psi-\psi_n=1$ on $U_n\cap\Omega$.
Then  $\varphi_n=0\vee(\psi-\psi_n)\wedge 1\in D(h_N)$, 
 because $h_N$ is  a Dirichlet form, $0\leq\varphi_n\leq1$, $\varphi_n=1$
 on $U_n\cap\Omega$ and 
$ \|\varphi_n\|_{D(h_N)}\leq \|\psi-\psi_n\|_{D(h_N)}\to0$
 as $n\to\infty$ again  by the Dirichlet property.
 Thus $\capp(B)=0$.

 \smallskip
 
 \noindent\ref{pcap2}$\Rightarrow$\ref{pcap1}$\;$
Assume   the balls $B(r)$ are bounded. Then  $(D(h_N)\cap L_\infty(\Omega))_c$  is a core of $h_N$
 by Proposition~\ref{pcap2.2}.
Therefore it  suffices to prove that each $\varphi\in (D(h_N)\cap L_\infty(\Omega))_c$ 
can be approximated in the $D(h_N)$-graph  norm by a sequence $\varphi_n\in (D(h_D)\cap L_\infty(\Omega))_c$.
If $A=(\supp\varphi)\cap\partial\Omega$ then $\capp(A)=0$ and  one may choose $\eta_n\in D(h_N)\cap L_\infty(\Omega)$ and open sets $U_n\subset \Ri^d$ such that $A\subset U_n$, $0\leq \eta_n\leq1$, $\eta_n=1$ on $U_n\cap\Omega$ and $\|\eta_n\|_{D(h_N)}\to0$ as $n\to\infty$.
Then set $\varphi_n=(\one_\Omega-\eta_n)\,\varphi$.
It follows that  $\varphi_n\in (D(h_D)\cap L_\infty(\Omega))_c$.
Moreover, by estimates similar to those used to prove Proposition~\ref{pcap2.10}
one deduces that $\|\varphi_n\|_{D(h_N)}\to0$ as $n\to\infty$.
Hence $H$ is Markov unique.\hfill$\Box$

\begin{cor}\label{ccap2.1}Assume $c_{kl}\in W^{1,\infty}_{\rm loc}(\overline\Omega)$ and that the balls $B(r)$ are bounded  for all $r>0$.
Then the following conditions are equivalent:
\begin{tabel}
\item\label{ccap1}
$\cc_A$ is satisfied for each bounded 
subset $A$ of $\overline \Omega$.
\item\label{ccap2}
$\capp(\partial\Omega)=0$.
\end{tabel}
\end{cor}
\proof\
It follows from Theorem~\ref{ntcap2.11} that Condition~\ref{ccap1} is equivalent to Markov uniqueness of $H$ and it follows from
Proposition~\ref{pcap} that Markov uniqueness of $H$ is equivalent to Condition~\ref{ccap2}. \hfill$\Box$ 

\bigskip

The proof of the corollary is indirect but if $\Omega$ is bounded then there is a simple direct proof
which shows that the two conditions of the corollary are complementary.
  Condition~\ref{ccap1}  is valid for bounded $\Omega$  if it is valid for $A=\overline\Omega$, i.e.\ the condition is equivalent to the existence of $\eta_n \in D(h_D)$ such
 that    $\lim_{n\to\infty} h_D(\eta_n) = 0$
and $\lim_{n\to\infty}\| \one_\Omega-\eta_n \|_2 = 0$.
Then, however, $\psi_n=\one_\Omega-\eta_n\in D(h_N)$, $\psi_n=1$ near $\partial\Omega$ and $\|\psi_n\|_{D(h_N)}\to0$ as $n\to\infty$.
Thus $\capp(\partial\Omega)=0$.
Conversely if $\capp(\partial\Omega)=0$ then there exist $\psi_n\in D(h_N)$ with $\psi_n=1$ near   $\partial\Omega$ such that  $\|\psi_n\|_{D(h_N)}\to0$ as $n\to\infty$.
Then setting $\eta_n=\one_\Omega-\psi_n$  one has $\eta_n \in D(h_D)$ and these functions  satisfy Condition~\ref{ccap1} of the corollary.

\section*{Acknowledgements} 
This work was was begun whilst the author was visiting the Mathematics
Department of Hokkaido University at the invitation of Akitaka Kishimoto.
It is a natural outcome of an ongoing collaboration with Adam Sikora to whom the author is indebted for many discussions and illuminating insights. 
The author also thanks El Maati Ouhabaz for some critical comments and Sasha Grigor'yan for   helpful correspondence.

\end{document}